# Discrete Euclidian Spaces: a starting point toward the discretization of mathematics

Ricardo Ramos Montero

*Departamento de Informática, Universidad de Oviedo*
Campus de Viesques, 33204 Gijón, Asturias, Spain

rramos@uniovi.es

Discrete Euclidian Spaces (DESs) are the beginning of a journey without return towards the discretization of mathematics. Important mathematical concepts—such as the idea of number or the systems of numeration, whose formal definition is currently independent of Euclidean spaces —have in the Isodimensional Discrete Mathematics (IDM) their roots in the DESs. This mathematics, which arises largely from the discretization of traditional mathematics, presents its foundations and concepts differently from the orthodox way, so at first glance it may seem that the IDM could be an exotic tool, or perhaps just "a simple curiosity." However, the IDM discrete approaches have a great theoretical repercussion on traditional mathematics.

## Introduction

Perhaps you have read, or you have wondered more than once, what makes a space of dimensionless points such as the Euclidean, locally modeling a universe where everything is dimensional. Would it not be more logical that the points of the mathematical models were also dimensional?

After spending several years developing a graphical system that works with *volumetric information*[1], the hierarchical spatial organization used with voxels [1] ended up being a three-dimensional discrete Euclidean space (3D-DES).

Since the theoretical study of the $n$D-DES generated high expectations, I have temporarily shelved the graphical system, working since then on the mathematical analysis of the $n$D-DES. As a result, I have established discrete mathematical foundations and concepts, giving rise to the *isodimensional discrete mathematics* (IDM), so called because the discrete points have the same dimension as the Euclidian space. As expected, the IDM stems largely from the discretization of traditional mathematics. For this reason, the most striking aspects of the IDM are those that differ from the established. In order to appreciate these conceptual differences we should begin by seeing what a discrete Euclidean space is in the IDM.

---

[1] The volumetric information is associated to the voxels, which are the three-dimensional equivalent of the famous pixels.

## Discrete Euclidean spaces

From the perspective of traditional mathematics, an $n$D-DES would be something like an "empty fractal grid", i.e., multiple grids organized using the views and concepts of fractals, but without any information registered on them. However, we can also obtain the definition of $n$D-DES in a simpler and clearer fashion, by taking advantage of the ideas and techniques of volumetric modelling, i.e., doing the discretization and the hierarchical arrangement of the usual Euclidean spaces.

### Hierarchical discretization

The hierarchical discretization of an $n$-dimensional continuous Euclidean space is carried out by dividing it into cells ($n$D-*points*), which are organized hierarchically according to their size. To perform the hierarchical discretization there are at least two methods: the bottom-up hierarchical arrangement and the top-down discretization.

The bottom-up hierarchical arrangement groups the $n$-dimensional discrete points according to a given criterion, starting with a *local DES*, i.e., a finite set of equal $n$D-points that, if glued, fill the whole Euclidean space[1]. Each group formed with the adjacent $n$D-points is a new $n$D-point, located at the next higher level. The same operation is repeated again and again with the newly formed $n$D-points, thus forming new growing discrete points, located at higher levels. At any scale level, the glued points cover the entire original Euclidean space, so that there is no excess or lack of space. The process ends when we have a single $n$D-point located at the root of the scale, with the same size as the original Euclidean space.

---

[1] A two-dimensional local DES is not graphically different from the grids.

On the other hand, the top-down discretization involves the subdivision of an $n$D-point (the original Euclidean space) into a finite number of cells, which become discrete $n$D-points at the level immediately below, setting in this way a new local DES at this level. In turn, these $n$D-points are divided again into smaller cells, all the same shape and size, becoming part of the next local DES, etc. Thus, as in the bottom-up methodology, different local DESs are generated in successive levels, hence the concept of scale is intrinsic to the $n$D-DES.

In short, an $n$D-DES can be understood as a spatial superposition of local DESs. For this reason, the study of the $n$D-DES is usually carried out on an $n$D-DES *breakdown*, which shows some of the local DESs from the scale root downward. In Figure 1 we see the first four local DESs, which correspond to a 2D-DES breakdown.

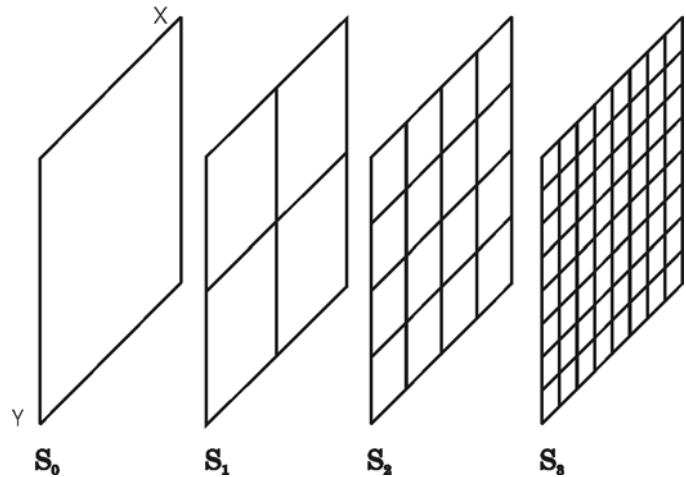

Figure 1: Four local DESs of a 2D-DES breakdown

We can observe that the 2D-points at any local DES are divided into four 2D-points in the next level, so that the 2D-DES of the example has a scale of order 4. The order of the scales is usually of the $2^k$ form, although there are no theoretical restrictions in this regard.

## Labeling the $n$D-DES

The definition of the $n$D-DES in any dimension would make little sense if it were not possible to recognize and to individually access each $n$D-point that forms a local DES, at a given level. To this end, once discretized the Euclidean space, the first thing we do is to properly label the resulting $n$D-DES.

The "natural" way of labeling the 2D-DES— or at least the one that offers major mathematical advantages —is shown in Figure 2.

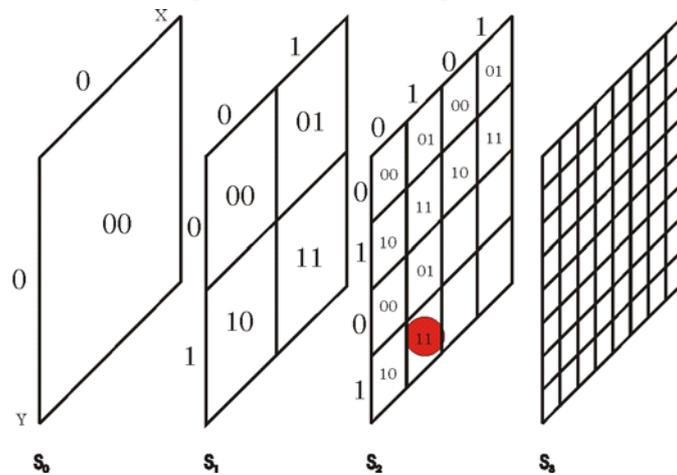

Figure 2: Labeling points on a 2D-DES breakdown

As we see, it is pure local labeling, that is, at any level of the scale the indexes are always the same. What advantages has the *local matrix* labeling? If we associate the prior indexes to a tree structure of the same order as the scale, we can better appreciate its virtues.

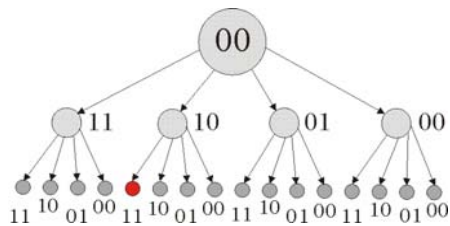

Figure 3: Labeling of the nodes in a tree of order 4

In this tree, each node at a level is clearly located by its *scalar index ($X_s$)*, that is, the list of local indexes of the nodes that form a given branch of the tree. For example, the scalar index of the point marked in Figure 2, and its corresponding node in Figure 3 is $X_s$ = (00 10 11).

Assuming that there is information associated with the $n$D-points of a local DES (e.g. that of an object), it could be recorded in the computer's memory using a tree structure like the previous one. Therefore, the information located in the $n$D-DES is easily recorded and located using trees of the same order as the spatial scale, if they are indexed as indicated. However, there is a problem.

No doubt, the scalar indexes are suitable to address trees such as the example, but we cannot say the same when they are indexing the local DESs, since they have a matrix arrangement. Thus, for a 2D-DES such as the one shown above, the use of Cartesian pairs (at each level) would be more appropriate. In Figure 2 the 2D-point marked in the $S_2$ level would normally be indexed by $(y, x) = (11, 01)_{radix\ 2} = (3, 1)_{radix\ 10}$. Is it possible to combine the scalar indexes with the Cartesian ones?

From Gargantini's work [2] we found a direct relationship between the two formats of indexing, as shown in Table 1.

| Level | $N_{t-1}$ | $N_{t-2}$ | ... | $N_0$ | |
|---|---|---|---|---|---|
| Y = | $y_{t-1}$ | $y_{t-2}$ | ... | $y_0$ | → Y Cartesian index |
| X = | $x_{t-1}$ | $x_{t-2}$ | ... | $x_0$ | → X Cartesian index |
| $X_s$ = | ↓<br>$(yx)_{t-1}$<br>matrix index at $N_{t-1}$ | ↓<br>$(yx)_{t-2}$<br>matrix index at $N_{t-2}$ | | ↓<br>$(yx)_0$<br>matrix index at $N_0$ | → 2D Scalar index |

Table 1: Relationship between scalar and Cartesian indexes

According to this relationship, the scalar index $X_s(yx)$ = (00 10 11) corresponds to the Cartesian pair $(y, x) = (011, 001)_2 = (3, 1)_{10}$, and vice versa. Is it always necessary to use binary indexes? Another example will show us that the answer is negative.

Let us suppose now that the dimensional points belong to a 3D-DES, with a scalar order of $4^3$, i.e., every time we go down a level in the scale, the 3D-points are subdivided into 64 3D-points. Then, in the tree of order 64, the scalar index node that corresponds to the 3D-point indexed by the Cartesian $(z, y, x)$ = (12, 15, 6)$_{10}$ = (030, 033, 012)$_4$ would be $X_s(zyx)$ = (000 331 032)$_4$ = (0 61 14)$_{10}$. As can be seen, changing from one format to another is very easy and fast, something of extreme importance in practical applications.

This important relationship between Cartesian and scalar indexes is not the only advantage offered by the $n$D-DES when working with them. For example, the partition of the indexes (typically Cartesians) allows working in a simple manner with the $n$D-points and their associated information defined on any scale level. Thus, in a 2D-DES with eight scale levels the points of the local DES in the last level *($S_7$)* would be indexed by the Cartesian pair *(y, x) = ($y_7$ … $y_0$, $x_7$ … $x_0$)*. But, what if we want to work at a higher level? Cartesian pair *(y, x) = ($y_7$, $x_7$)* would work without problems at the $S_1$ level, with the *$(b^2)^1$* 2D-points of this level, being *$b^2$* the order of the scale. If you then wish to operate in $S_3$, the Cartesian *(y, x) = ($y_7 y_6 y_5$, $x_7 x_6 x_5$)* provides access to the *$(b^2)^3$* 2D-points of the level. Therefore, a simple shift of indexes is enough to get the information associated with the $n$D-points on any local DES[1].

Although there are many advantages in working with information associated with the $n$D-DES, our goal is not to show the practical skills of these spaces, but to see why they are the starting point toward the discretization of traditional mathematics. To this end, first we will see how the information is associated with the dimensional points of the $n$D-DES.

## Information associated with the dimensional points

The $n$-dimensional points of a local DES can be of three types, depending on the associated information[2]. Thus, an $n$D-point will be [*absolute*] *homogeneous*, if all the information that we find at the sublevels of the scale rooted at that point is of the same type. For instance, in a scale of order $b^2$, assuming that an $n$D-point at the $S_k$ level is red, if the $(b^2)^1$ points at $S_{k+1}$ are red too, and if the $(b^2)^2$ points at $S_{k+2}$ are also red, and so on, then we can say that the $n$D-point at the $S_k$ level is absolute homogeneous. On the contrary, if we find one or more empty $n$D-points (i.e., without associated information) at any sublevel of the scale or with different information, then our $n$D-point at the $S_k$ level will be

---

[1] This modus operandi might be called *true floating-point calculation,* since the dot in the standard floating-point does not move at all in the scale, although apparently it does. Unfortunately, the *true floating-point calculation* nowadays is not as efficient as it could be due to the fact it is not implemented in the hardware of computers.

[2] For practical purposes there are three types of points, but from a theoretical point of view they are only two.

*relative homogeneous*, if it is located at the terminal level (the lowest level of the scale). Otherwise it would be heterogeneous[1].

These concepts, taken to a tree structure of order $b^2$ suggest that homogeneous $n$D-points (absolute or relative) will always be represented by terminal nodes (or leaf nodes), while nodes at higher levels, up to the tree root will be heterogeneous. Note that the presence of relative homogeneous information is absolutely necessary. Otherwise, it would not be feasible to record most of the information we want to.

These ideas on how to attach and record information on the $n$D-DES are crucial to understand the concept of number in the isodimensional discrete mathematics, as we shall now see.

## The concept of number in the IDM

Let us review the ideas we have seen about the $n$D-DES and the associated information, but this time focusing on the 1D-DES, and using a new terminology. There is no better way to do this than to begin showing the breakdown of a 1D-DES of order 4, and the appropriate labeling (Figure 4). In this dimension we will often use the term *scalar perpendicular* or *numerical sequence*, to refer to the *scalar index*.

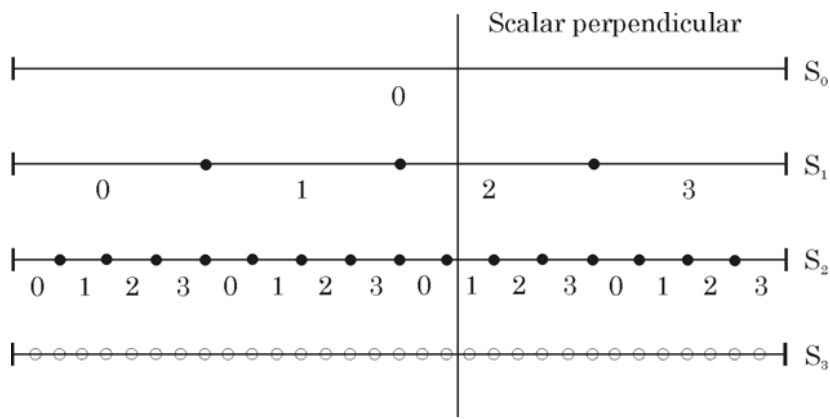

Figure 4: Breakdown and labeling of a 1D-DES

---

[1] Terminal 3D-points with associated information about matter are always relative homogeneous (at least in the known scalar levels), because as we move down the scale, we find different information about the elements or constituent particles (e.g., water, hydrogen and oxygen atoms, electrons, neutrons, protons, etc.). Therefore, a terminal 3D-point with "water" as associated information would be relative homogeneous.

Assuming the scale is open at the bottom, i.e., with an unlimited number of subscalar levels, the scalar perpendicular in Figure 4 provides the numerical sequence $X_s$ = 021... Is this scalar index a number? Not yet, because for a numerical sequence to become a number, *it must be evaluated*, i.e., we must assign to it a numerical value (numerical information). To this end, we have to choose the *evaluation level*. In doing this, the numerical sequence will be automatically evaluated and converted into a number. Which is the value associated with the sequence? The numerical value of a number is the quantity of 1D-points that there are from the origin of the local DES of evaluation, to the 1D-point crossed by the scalar perpendicular, exclusive. For example, if we choose $S_2$ as the level of evaluation, we have the number $X = (021)_4 = (9)_{10}$, as we can verify directly above. Is it a coincidence? No, the remarkable fact is that the Cartesian and scalar indexes are the same in 1D-DES[1]. After evaluating a numerical sequence it becomes a number but, what kind of number? As we have seen previously, the terminal $n$D-point with information associated can be absolute or relative homogeneous. Is it the same when we assign numerical information? Yes, it is, but using different terminology.

So, if there is no numerical information in the subscale rooted at the 1D-point crossed by the scalar perpendicular at the evaluation level, the 1D-point is absolute homogeneous. In the usual terminology this means that the associated value is an integer (or natural), and the numerical sequence is terminal. For instance, the sequence 02011,0... is terminal, because there is no numerical information from the evaluation level (indicated by a comma). For the contrary, 02011,5738... is a non-terminal sequence, since the 1D-point at the level of evaluation is relative homogeneous, i.e., the subscale with origin in that 1D-point has more numerical information. The concepts of terminal and non-terminal, or evaluated and non-evaluated numerical sequences are of crucial importance in the theoretical foundations of the IDM [4].

As can be seen, it is difficult to imagine a more natural way to establish the concept of number, with the added advantage that the numbers stem from the space itself, i.e., the place where they work perfectly, measuring and evaluating all kinds of objects and events. In addition, the scalar vision of numbers pro-vided by the discrete Euclidean spaces[1] allows another perspective of the decimal point, since it no longer determines not only the value of the numbers, but also the level of observation or appreciation of the calculations (another way of looking at the numerical value). For example, the sum of the terminal sequences 001.234,0... and 011.283,0... is carried out without doubt at $S_5$ level, resulting 012.517,0...[2]. The dot just indicates that we see (or appreciate) the calculation from $S_2$[3]. These ideas, taken to the field of *analysis* extend the scope of *differential and integral calculus* [4].

---

[1] This is why the conversion from one format to another (shown in Table 1) is so simple.

## Positional numeral systems

The 1D-DES not only allows a natural establishment of numbers, but also all concepts related to them also arise naturally from this discrete space, as for instance the positional numeral systems or the justification of the existence of zero.

Indeed, each digit of a numerical sequence indicates the number of complete 1D-points that the $S_i$ level provides to the value of the sequence. Obviously, this contribution to the value of the sequences is given in $S_i$, so this local value has to be adapted to the evaluation level *($S_e$)* in order to obtain a scalar homogeneous result, a process that depends on the order of the scale. For example, being $X_s$ = 021 the numerical sequence (Figure 4), we see that the 1D-point at $S_0$ (root) does not provide any complete 1D-points to the value of the sequence, since it is crossed by the scalar perpendicular. In contrast, the $S_1$ level contributes with two 1D-points (0 and 1), which are converted into $2 \cdot (4)^1$ 1D-points at the evaluation level *($S_2$)*, being 4 the scalar order. Finally, the evaluation level itself gives an extra 1D-point (the 0), since the 1D-point indexed by 1 is split by the perpendicular. In short, the value of the sequence evaluated at $S_2$ (021) is given by $0 \cdot 4^2 + 2 \cdot 4^1 + 1 \cdot 4^0 = 9$. Note that whenever the scalar perpendicular goes through a 1D-point indexed by 0, the contribution to the value of the numerical sequence is null.

---

[1] In many cases, it is more interesting to see numbers as scalar sequences or indexes than consider them just by their associated value.

[2] In this example, it is clear why the floating-point is actually a "fixed-point" calculation, since the scalar homogeneity between operands and the result is necessary to calculate consistently.

[3] In other words, the decimal point makes a partition in the Cartesian and scalar indexes.

What if we eliminate the 0, labeling the scale after that from 1 to *a*, being 10 the value of the index *a*? The result is also a positional numeral system ([3], [4]), which is as valid to count as the decimal system, but it is not appropriate to indexing the 1D-DES. Furthermore, the result of *(t − t)* is not represented in the system. In short, the zero, with a null contribution to the value of the evaluated numerical sequences is absolutely essential in the *optimal* positional numeral systems. How many optimal positional numeral systems are there? Merely changing the order of the scale is enough to obtain a different positional numeral system[1]. What if we change to higher dimensions?

If we go to a 2D-DES, the positional numeral systems in this space hardly differ conceptually from what we have just seen, although the numerical sequences (i.e., the scalar indexes) are of the form *ab cd ef...* as shown above. Therefore, each matrix index *(yx)* is now a figure of a number expressed on a two-dimensional positional numeral system. The process of homogenization at the level of evaluation is the same. Thus, after evaluating the scalar index that is shown in Figure 2, we get the two-dimensional number $X(yx) = (00\ 10\ 11,)_2$. Its numerical value— i.e., the total of 2D-points that are in the local DES at $S_2$, following the indexing order —is given by $[00 \cdot (10^2)^2 + 10 \cdot (10^2)^1 + 11 \cdot (10^2)^0]_2$, that expressed in the decimal numeral system becomes $0 + 2 \cdot 4 + 3 \cdot 1 = 11$. Interestingly enough, the two and three-dimensional positional numeral systems are widely used in mathematics, for instance working with units of area or volume, or working out square roots by hand[2]. However, textbooks usually do not explicitly specify them.

Before concluding this section it is worth comparing this way of seeing and understanding the concept of number, with the one that presents the *axiomatic set theory*. To do this, we will simply analyze the zero in both contexts.

As you know, the axiomatic set theory associates the empty set to zero, that is to say it is given a null value. This approach has at least two drawbacks:

1) The zero's null value is established by axiomatic definition.

---

[1] Note that in one dimension the order of the scale and the radix are the same. In higher dimensions this is not so.

[2] This algorithm requires the figures of the radicand grouped in pairs from the dot.

2) The explicit assignment of the null value has contributed significantly to the inhibition of the zero in its indexing role[1].

Note that the IDM in principle does not define the numerical value of zero. It simply confirms that the contribution of a local DES to the value of the evaluated numerical sequence is null, if the 1D-point crossed by the scalar perpendicular is indexed by 0. In any case, does the figure 0 have a null value? Counting the complete 1D-points that the $S_0$ level provides to the evaluated numerical sequence $X = (0,)$, it is obvious that the value of $X$ is null. However, the evaluated sequence $(0,)$ and the non-evaluated sequence $(0)$ theoretically are not the same concept. The first is a number, and the second is an index (label), so that it is incorrect to say that the figure 0 has a null value, but it is mathematically acceptable because the *number zero* $(0,)$ has. Then, this is just a small sample of the important and subtle differences between both types of maths.

## The discrete line

Traditional mathematics ensures that between two different points on the real line, regardless of how close they are, there are infinitely many points. The IDM can also make a similar claim in the scope of the *discrete line*, but taking into account the dimensionality of the points.

First, by "discrete line" we shall understand a local DES consisting of an unlimited number of 1D-points. Therefore, we have to assume that it is located at the $S_\infty$ level of a 1D-DES with an endless scale of order $b$, i.e., without a subscalar limit.

The second issue is that, between two consecutive 1D-points of any local DES at $S_k$, and of course in the discrete line too, there are no dimensional points. However, between two alternative 1D-points there are $b^0$ 1D-point at $S_k$, $b^1$ at $S_{k+1}$, $b^2$ at $S_{k+2}$, and $b^\infty$ 1D-points at the $S_{k+\infty}$ level. Therefore, between two alternatives 1D-points of the discrete line at $S_k$ there is an unlimited number of 1D-points defined at the $S_\infty$ level.

Although the concepts of discrete and real line are not the same, in daily practice there is no difference between them, since numbers with $t$ decimal

---

[1] One of the traditional faults in mathematics at all times has been to start indexing from 1.

figures in the interval [0, 1] are the same in both lines, and with the same numerical value[1].

## Numeric systems

So far, we have only seen a few examples of how the IDM addresses some mathematical concepts, but I think they are enough to understand that it is not possible to discretize mathematics keeping the current foundations.

As expected, the adoption of new foundations in the IDM affects the majority of theories and fields of mathematics. As a small sample to show that this is true, we will see now that something as basic as the definition of the usual *number systems* does not fit in the IDM.

The orthodox mathematics classifies the different types of numbers in 5 basic sets: *naturals (N), integers (Z), rational numbers (Q), real numbers (R) and complex numbers (C)*. Each of these sets includes or comprises the former, so that they differ only in a subset of numbers that is not present in previous sets. For example, the set of real numbers includes all the previous sets, plus the irrational numbers.

If we try to organize these number systems in the IDM applying the same criteria, failure is assured, except in the case of naturals that match in both contexts. Integers do not fit completely, since the way the IDM sees the negative numbers differs from that of traditional maths. Finally, on establishing the rational numbers the respective views are irreconcilable. I will consider this case, because it is easy to explain.

As you know, rational numbers in mathematics are those that can be expressed as a fraction of integers, i.e., if $a$, $b$ are integers, then the result of $a/b$ is a rational number. This means that all numerical sequences of rational numbers are periodic or terminal sequences, since it is in the nature of the process of division to generate this kind of sequences.

Accordingly, the sequence 0.9... should be a rational number, because it is periodic. However, a generating fraction for this sequence does not exist, therefore it cannot be a rational number. In mathematics this issue is not a problem, since it can be proved that 0.9... = 1. In contrast, the IDM cannot assert such a thing, because a simple glance at the 1D-EDE breakdown of order 10 is enough to convince us that 0.9... ≠ 1. Thus, we have a case of pure-periodic sequence and an unlimited amount of mixed-periodic sequences, which cannot be classified as rational numbers, and certainly they are not irrationals.

---

[1] The evaluated numerical sequences have by definition a limited number of digits. This implies that any numerical comparison between both maths has to be made among truncated numbers. The theoretical implications of working with endless figures go beyond the scope of this article [4].

## Conclusions

The IDM could be seen as an odd or exotic way of setting the foundations and mathematical concepts. However, these approaches must be taken into account, because they have an enormous theoretical impact on mathematics. Why? Simply because the IDM provides another point of view about concepts firmly established for centuries, and what is more important, IDM brings another rigorous way of doing mathematics. Besides, IDM opens new paths in basic mathematical research on topics that many have given as settled a long time ago.